# GAUSS-MARKOV PROCESSES ON HILBERT SPACES

BEN GOLDYS, SZYMON PESZAT, AND JERZY ZABCZYK

Abstract. K. Itô characterised in [13] zero-mean stationary Gauss Markov-processes evolving on a class of infinite-dimensional spaces. In this work we extend the work of Itô in the case of Hilbert spaces: Gauss-Markov families that are time-homogenous are identified as solutions to linear stochastic differential equations with singular coefficients. Choosing an appropriate locally convex topology on the space of weakly sequentially continuous functions we also characterize the transition semigroup, the generator and its core thus providing an infinite-dimensional extension of the classical result of Courrège [3] in the case of Gauss-Markov semigroups.

Contents

## 1. Introduction

The aim of the paper is to derive several characterizations of Gauss-Markov processes on infinite dimensional Hilbert spaces.

1991 *Mathematics Subject Classification.* 60G15, 60H15, 60J99.

The work of Ben Goldys was partially supported by the ARC Discovery Grant DP120101886. Part of this work was prepared during his visit to the Institute of Mathematics of the Polish Academy of Sciences. Ben Goldys gratefully acknowledges excellent woking conditions and stimulating atmosphere of the Institute.

The work of Szymon Peszat and Jerzy Zabczyk was supported by Polish Ministry of Science and Higher Education Grant N N201 419039 "Stochastic Equations in Infinite Dimensional Spaces".

Let $H$ be a real separable Hilbert space and let $\{\mu(t,x)\colon x\in H,\, t>0\}$ be a family of *Gaussian transition kernels* on $H$, that is

$$\mu(t,x)=\mathscr{N}\left(m(t,x),Q(t,x)\right), \tag{1.1}$$

and for every Borel set $B\subset H$,

$$\mu(s+t,x)(B)=\int_H \mu(s,x)(\mathrm{d}y)\mu(t,y)(B),\quad s,t\geqslant 0,\ x\in H. \tag{1.2}$$

It is well known, see for example Theorem 7.4 in [14], that there exist a measurable space $(\Omega,\mathscr{F})$ endowed with a filtration $(\mathscr{F}_t)$, an $H$-valued and $(\mathscr{F}_t)$-adapted stochastic process $Z$ defined on $H$, and a family of measures $\{\mathbb{P}^x\colon x\in H\}$ on $\Omega$ such that for every $x\in H$ the process $Z$ is Markov on the probability space $(\Omega,\mathscr{F},(\mathscr{F}_t),\mathbb{P}^x)$, its transition kernel is $\mu(t,x)$ and $\mathbb{P}^x\left(Z(0)=x\right)=1$. Let us recall that the system $(\Omega,\mathscr{F},(\mathscr{F}_t),Z,\{\mathbb{P}^x\colon x\in H\})$ is called a *Markov family*.

Our aim in this paper is to derive three explicit characterizations of such processes. Our Theorem 2.2 describes the structure of the transition semigroup, Theorem 2.3 provides a construction of a stochastic equation which the process satisfies and Theorem 2.4 characterizes the generator of the transitions semigroup. In particular, we will prove that there exist

- a $C_0$-semigroup $(L(t))$ (with the generator $A$) acting on $H$,
- a vector $b_H=b_H(\lambda)\in H$ defined for a certain $\lambda>0$,
- a selfadjoint operator $Q\geqslant 0$ in $H$,

such that for any $h\in\operatorname{dom}(A^\star)$ the mean value and the covariance operator of $Z(t,x)$ satisfy for certain $\lambda>0$,

$$\langle m(t,x),h\rangle=\langle L(t)x,h\rangle+\int_0^t \langle L(s)b_H,(\lambda-A^\star)\,h\rangle\,\mathrm{d}s,\quad t\geqslant 0,\, x\in H,$$

and

$$\langle Q(t,x)h,h\rangle=\int_0^t \left|Q^{1/2}L^\star(s)h\right|^2\mathrm{d}s,\quad t\geqslant 0,\, x\in H.$$

Given this representation of $m$ and $Q$ we derive the second and the third characterizations. The second one is an extension of the work of K. Itô who obtained in [13] a similar representation for a stationary Gauss-Markov process with zero mean. Since we do not assume stationarity and we work with a Gauss-Markov family of processes instead of a single process, the approach of [13] can not be adopted and a different argument is required to derive representations for $m(t,x)$ and $Q(t,x)$. Given these representations, we prove that there exists an $(\mathscr{F}_t)$-adapted standard cylindrical Wiener process on $H$ such that

$$\begin{cases}\mathrm{d}Z=(AZ+b_V)\,\mathrm{d}t+Q^{1/2}\mathrm{d}W,\\ Z(0)=x\in H,\end{cases} \tag{1.3}$$

where $b_V=(\lambda-A)b_H$ is uniquely defined in a larger Hilbert space $V\supset H$ to which the semigroup $(L(t))$ can be extended.

Our third, purely analytic result is in the spirit of a theorem of Ph. Courrège [3]. Namely, Courrège obtained a representation of generator for an arbitrary Markov semigroup, that is strongly continuous in $C_0\left(\mathbb{R}^d\right)$[1] and such that $C_c^\infty\left(\mathbb{R}^d\right)$ is a core for its generator. It is well known, see for example [5] that the space $C_b^{norm}(H)$ of bounded continuous functions on $H$ endowed with the supremum norm is not appropriate for the analysis of Markovian semigroups on Hilbert spaces. In order to obtain a counterpart of the Courrège theorem we need to introduce two locally convex topologies. The first one is the bounded weak topology on $H$ as defined in [10]. It was demonstrated in [17] and recently in [2] that, in fact, it is a natural concept for Markovian transition semigroups and particularly useful for proving the existence of an invariant measure. We will use the notation $H_{bw}$ for the space $H$ endowed with the bounded weak topology, $C_b\left(H_{bw}\right)$ for the space of bounded continuous functions on $H_{bw}$, and $C_b^{norm}\left(H_{bw}\right)$ for the space $C_b\left(H_{bw}\right)$ endowed with the norm topology. It was shown in [17] that $C_b\left(H_{bw}\right)$ is precisely the space of bounded functions that are weakly sequentially continuous. Another locally convex topology, known as the strict topology, should be introduced on the space $C_b\left(H_{bw}\right)$. The strict topology can be defined as the strongest topology on $C_b\left(H_{bw}\right)$ which on norm bounded sets of $C_b^{norm}\left(H_{bw}\right)$ coincides with the topology of uniform convergence on compacts of $H_{bw}$. For alternative definitions of this topology and its applications to the theory of Markov semigroups, see [12] and references therein.

Using the representations for $m(t,x)$ and $Q(t,x)$ we show that the transition semigroup $P_t\phi(x) = \mathbb{E}\phi(Z(t,x))$ of the process $Z$ is strongly continuous in $C_b^{strict}\left(H_{bw}\right)$. Moreover, we will identify the core of the generator $L$ in $C_b^{strict}\left(H_{bw}\right)$ and will derive an explicit formula for the generator acting on functions from the core. In the case when $A$ has bounded inverse this formula takes a simple form

$$L\phi(x) = \frac{1}{2}\operatorname{Tr}\left(QD^2\phi(x)\right) + \left\langle x + b_H, A^\star D\phi(x)\right\rangle, \quad x \in H.$$

For a general version of this formula see Section 2. These results provide a version of the Courrège theorem for a Gauss-Markov semigroup in a Hilbert space.

Although the extension concerns a rather narrow class of transition semigroups we believe that it does suggest how a general infinite dimensional version of the Courrège theorem should look like.

If the state space $H$ is one-dimensional then results of our paper are part of the mathematical folklore. More generally, it is not very difficult to obtain our result if $\dim(H) < \infty$. Let us note here that in finite dimensional spaces a similar problem is considered in the framework of affine processes, see for example [9]. In the theory of affine models a much wider class of processes is considered but the linearity of the function $x \mapsto m(t,x)$ and the fact that $x \mapsto Q(t,x)$ is constant are assumed from the very beginning while we do derive these properties as a result of a careful analysis.

[1]We use standard notations $C_0\left(\mathbb{R}^d\right)$ for the space of continuous functions vanishing at infinity, and $C_c^\infty\left(\mathbb{R}^d\right)$ for the space of infinitely differentiable functions with compact supports.

In infinite dimensions the only result in this directions was obtained by Itô in [13]. Let us also recall a related problem to characterise a family of measures that satisfy the so-called skew convolution equation [19]. In that paper it is assumed again from the very beginning that the expectation $m(t,x)$ is linear in $x$ and the covariance $Q(t,x)$ constant in $x$.

Main results are formulated in Section 2. The proofs are presented in the following sections. In the final Section 6 we discuss an application of our results to models with boundary noise.

## 2. Main results

We give now precise formulation of our main theorems. The proofs will be postponed to the following sections.

By definition we have

$$m(t,x)=\mathbb{E}^x Z(t),\quad x\in H,\, t\geqslant 0,$$

and for any $h,k\in H$

$$\langle Q(t,x)h,k\rangle=\mathbb{E}^x\left\langle Z(t)-m(t,x),h\right\rangle\left\langle Z(t)-m(t,x),k\right\rangle,\quad x\in H,\ t\geqslant 0.$$

Let $\phi: H\to\mathbb{R}$ be continuous and such that $|\phi(x)|\leqslant C\left(1+|x|^2\right)$. Then for any $s,t\geqslant 0$ the Markov property yields

$$\mathbb{E}^y\left(\phi\left(Z(t+s)\right)\middle| Z_s=x\right)=\mathbb{E}^x\phi\left(Z(t)\right)\quad \mathbb{P}^y-a.s.,\quad x,y\in H.$$

The following hypothesis is a standard assumption of the rest of the paper.

**Hypothesis 2.1.** *(1) For every $x\in H$ and $h,k\in H$ the functions*

$$t\mapsto\langle m(t,x),h^\star\rangle\quad\text{and}\quad t\mapsto\langle Q(t,x)h^\star,k^\star\rangle$$

*are continuous.*
*(2) For any $t\geqslant 0$ and $h,k\in H$ the functions*

$$x\to\langle m(t,x),h\rangle\quad\text{and}\quad x\to\langle Q(t,x)h,k\rangle$$

*are continuous on $H$.*
*(3) For any $t>0$ and any $x\in H$,*

$$\overline{\operatorname{Im} Q(t,x)}=H.$$

We note that for every $x\in H$, Hypothesis 2.1 yields

$$\lim_{t\downarrow 0}\langle m(t,x),h\rangle=\langle x,h\rangle,\quad h\in H,$$

and

$$\lim_{t\downarrow 0}\langle Q(t,x)h,k\rangle=0,\quad h,k\in H.$$

In order to formulate our results we need some preparations. Let $L=(L(t))$ be any $C_0$-semigroup on $H$ and let $A$ denote its generator. Then for $\lambda>0$ big enough $(\lambda-A)$ is boundedly invertible on $H$ and a new norm on $H$ can be defined by the formula

$$|x|_{-1}=\left|(\lambda-A)^{-1}x\right|,\quad x\in H.$$

The space $V$ is defined as a completion of $H$ with respect to the norm $|\cdot|_{-1}$. Clearly, the topological space $V$ does not depend on $\lambda > 0$. The semigroup $(L(t))$ extends to a $C_0$-semigroup $(L_V(t))$ on $V$ with the generator $A_V$ whose domain is equal to $H$.

We can formulate now our first main result. It gives a precise description of the functions $m$ and $Q$ and will also play the crucial role in the proofs of Theorems 2.3 and 2.4 below. The proof of Theorem 2.2 is postponed to Section 3.

**Theorem 2.2.** *There exists a strongly continuous semigroup $(L(t))$ with the generator $A$ on $H$, a vector $b_H \in H$ and a selfadjoint operator $Q \geqslant 0$ in $H$ such that:*

(1) *We have*

$$\operatorname{dom}(A^\star) \subset \operatorname{dom}\left(Q^{1/2}\right). \tag{2.1}$$

(2) *For any $x \in H$,*

$$m(t,x) = L(t)x + \int_0^t L_V(s) b_V \, \mathrm{d}s, \quad t \geqslant 0, \tag{2.2}$$

*where $b_V = A_V b_H \in V$. In particular, for every $h \in \operatorname{dom}(A^\star)$,*

$$\langle m(t,x), h\rangle = \langle L(t)x, h\rangle + \int_0^t \langle L(s) b_H, (\lambda - A^\star)\, h\rangle \, \mathrm{d}s, \quad t \geqslant 0,\ x \in H. \tag{2.3}$$

(3) *The covariance operator $Q(t) = Q(t,x)$ is independent of $x \in H$,*

$$\operatorname{Im}(L^\star(t)) \subset \operatorname{Dom}\left(Q^{1/2}\right), \quad t > 0,$$

*and*

$$Q(t) = \int_0^t \left(Q^{1/2} L^\star(s)\right)^\star \left(Q^{1/2} L^\star(s)\right) \mathrm{d}s. \tag{2.4}$$

*In particular,*

$$\langle Q(t)h, h\rangle = \int_0^t \left|Q^{1/2} L^\star(s) h\right|^2 \mathrm{d}s, \quad t \geqslant 0,\ x \in H, \tag{2.5}$$

*and*

$$\int_0^T \left\|Q^{1/2} L^\star(t)\right\|_{HS}^2 \mathrm{d}t < \infty.$$

(4) *For any $x \in H$, any $h, k \in H$ and $0 \leqslant s \leqslant t$,*

$$\langle Q(s,t,x)h, k\rangle = \mathbb{E}^x \langle Z(s) - m(s,x), h\rangle \langle Z(t) - m(t,x), k\rangle = \langle L(t-s) Q(s) h, k\rangle, \tag{2.6}$$

*hence the operator $Q(s,t) = Q(s,t,x)$ is independent of $x \in H$.*

Now we can formulate our second main result. Its proof is postponed to Section 4.

**Theorem 2.3.** *Assume that Hypothesis 2.1 holds. Then there exist a strongly continuous semigroup of operators $(L(t))$ on the space $H$ with its associated space $V$, an element*

*$b_V \in V$, a selfadjoint and non-negative operator $Q$ in $H$, and there exists a standard cylindrical Wiener process $W$ on $H$ such that for every $x \in H$ we have*

$$Z(t) = L(t)x + \int_0^t L_V(t-s)b_V \mathrm{d}s + \int_0^t L(t-s)Q^{1/2}\mathrm{d}W(s), \quad t \geqslant 0,\ \mathbb{P}^x - a.s. \tag{2.7}$$

Our final result is an extension of the Courrège theorem. We introduce first certain convenient locally convex topologies.

In the terminology of [10] the bounded weak topology $\tau^{bw}$ on $H$ is the strongest locally convex topology on $H$ that coincides with the weak topology on every ball. We will denote by $H_{bw}$ the space $H$ endowed with the topology $\tau^{bw}$. It has been proved in [17] that $\phi \in C_b(H_{bw})$ if and only if $\phi$ is bounded and weakly sequentially continuous on $H$.

Let $\tau^n$ and $\tau^c$ denote, respectively, the norm topology and the topology of the uniform convergence on compacts on the space $C_b(H_{bw})$. The strict topology $\tau^s$ is defined as the strongest locally convex topology on $C_b(H_{bw})$ that coincides with the topology $\tau^c$ on compacts. We will use a short notation $C_b$ for the space $C_b(H_{bw})$ endowed with the topology $\tau^s$.

We will say that $\phi \in \mathscr{F}C_b^\infty(A^\star)$ if there exist $n \geqslant 1$, $f \in C_b^\infty(\mathbb{R}^n)$ and $h_1, \ldots, h_n \in \operatorname{dom}(A^\star)$ such that

$$\phi(x) = f(\langle x, h_1\rangle, \ldots, \langle x, h_n\rangle)$$

and

$$\sup_{x\in H} |\langle x, A^\star D\phi(x)\rangle| < \infty.$$

It can be easily shown that $\mathscr{F}C_b^\infty(A^\star)$ is dense in $C_b^{strict}(H_{bw})$. The proof of the result below is postponed to Section 5.

**Theorem 2.4.** *The semigroup $(P_t)$ is strongly continuous in $C_b^{strict}(H_{bw})$. Let $L$ denote the generator of the semigroup $(P_t)$ in $C_b^{strict}(H_{bw})$. Then the space $\mathscr{F}C_b^\infty(A^\star)$ is a core for $L$ and for every $\phi \in \mathscr{F}C_b^\infty(A^\star)$*

$$L\phi(x) = \frac{1}{2}\mathrm{Tr}\left(QD^2\phi(x)\right) + \langle x, A^\star D\phi(x)\rangle + \langle b_H, (\lambda - A^\star) D\phi(x)\rangle. \tag{2.8}$$

*Remark* 2.5. Let us note that $b_H$ depends on $\lambda$ but the formula (2.8) uniquely defines $L\phi$ since $b_V = (\lambda - A)b_H$ does not depend on $\lambda$.

## 3. Proof of Theorem 2.2

We recall first a basic fact about the conditional Gaussian measures.

3.1. **Conditional Gaussian measures.** Let $H_1$ and $H_2$ be two real separable Hilbert spaces and let $(X, Y)$ be a Gaussian vector defined on a probability space $(\Omega, \mathscr{F}, \mathbb{P})$ and taking values in $H_1 \times H_2$. Let

$$m_X = \mathbb{E}X \quad \text{and} \quad m_Y = \mathbb{E}Y.$$

The covariance operator $C_X$ of $X$ is determined by the equation

$$\mathbb{E}\langle X - m_X, h\rangle \langle X - m_X, k\rangle = \langle C_X h, k\rangle, \quad h, k \in H_1, \tag{3.1}$$

and a similar condition determines the covariance $C_Y$ of $Y$. The covariance operator $C_{XY}\colon H_1 \to H_2$ is defined by the condition

$$\langle C_{XY}h,k\rangle = \mathbb{E}\,\langle X-m_X,h\rangle\,\langle Y-m_Y,k\rangle\,,\quad h\in H_1, k\in H_2,$$

and then $C^*_{XY} = C_{YX}$. For a linear closable operator $G$ on $H$ the closure of $G$ will be denoted by $\overline{G}$. For a symmetric and compact operator $K\colon H\to H$ we will denote by $K^{-1}$ its pseudo inverse. For the convenience of the reader we present the following known result (see. e.g. [16]).

**Theorem 3.1.** *The following holds.*
*(a) We have*

$$\mathrm{Im}\,(C_{YX}) \subset \mathrm{Im}\left(C_X^{1/2}\right), \tag{3.2}$$

*and the operator $K := C_X^{-1/2}C_{YX}$ is of Hilbert–Schmidt type on $H$ and $K^* = \overline{C_{XY}C_X^{-1/2}}$.*
*(b) We have*

$$\mathbb{E}\,(Y|X) = m_Y + K^*C_X^{-1/2}\,(X-m_X)\,,\quad \mathbb{P}-a.s.$$

*(c) The conditional distribution of $Y$ given $X$ is Gaussian $\mathscr{N}\left(\mathbb{E}\,(Y|X)\,,C_{Y|X}\right)$, where*

$$C_{Y|X} := C_Y - K^*K.$$

*Moreover, the random variables $K^*C_X^{-1/2}X$ and $\left(Y-K^*C_X^{-1/2}X\right)$ are independent.*

3.2. **Proof of Theorem 2.2.** Let $Q(s,t,y)$ denote the covariance operator,

$$\langle Q(s,t,y)h,k\rangle = \mathbb{E}^y\,\langle Z(s)-m(s,y),h\rangle\,\langle Z(t)-m(t,y),k\rangle\,,\quad h,k\in H.$$

By (3.2) for any $y\in H$ and any $0<s\leqslant t$,

$$K(s,t,y) := Q^{-1/2}(s,y)Q^\star(s,t,y)$$

is a well defined Hilbert–Schmidt operator on $H$. The proof is divided into a sequence of lemmas.

**Lemma 3.2.** *(1) For any $s,t$ such that $0\leqslant s\leqslant t$ the operator-valued mappings $x\mapsto Q(s,t,x)$ and $x\mapsto K^\star(s,t,x)K(s,t,x)$ are constant in $x\in H$ and for $Q(t)=Q(t,t,0)$ and $K(s,t)=K(s,t,0)$ we have*

$$Q(t-s) = Q(t) - K^\star(s,t)K(s,t), \tag{3.3}$$

*(2) For any $s,t$ such that $0\leqslant s\leqslant t$ the operator $L(t,s) := K^\star(s,t)Q^{-1/2}(s)$ with the domain $\mathrm{Im}\left(Q^{1/2}(s)\right)$ has a unique extension to a bounded linear operator $L(t,s)\colon H\to H$ and for any $x,y\in H$,*

$$m(t-s,x) = m(t,y) - L(t,s)m(s,y) + L(t,s)x. \tag{3.4}$$

*Moreover,*

$$m(t-s,m(s,y)) = m(t,y),\quad 0\leqslant s\leqslant t,\; y\in H. \tag{3.5}$$

*Proof.* Invoking Theorem 3.1 we find that for any $0 < s < t$ and $y \in H$ the following equality holds $\mathbb{P}^y$-a.s.

$$\begin{aligned} &\mathbb{E}^y\left((Z(t)-\mathbb{E}^y(Z(t)|Z(s)))\otimes(Z(t)-\mathbb{E}^y(Z(t)|Z(s)))|\,Z(s)\right)\\ &\quad = Q(t,y)-K^\star(s,t,y)K(s,t,y). \end{aligned} \tag{3.6}$$

Applying the Markov property to (3.6) we obtain for $0 < s < t$,

$$Q(t-s,x) = Q(t,y) - K^\star(s,t,y)K(s,t,y) \quad \text{for} \quad \mu(s,y)-a.e.\ x. \tag{3.7}$$

For any $h, k \in H$,

$$\langle Q(t-s,x)h,k\rangle = \langle Q(t,y)h,k\rangle - \langle K^\star(s,t,y)K(s,t,y)h,k\rangle \quad \text{for} \quad \mu(s,y)-a.e.\ x. \tag{3.8}$$

Therefore the function $x \mapsto \langle Q(t-s,x)h,k\rangle$ is constant on a dense set and by Hypothesis 2.1 is continuous on $H$. It follows that equality (3.8) holds for every $x \in H$. Finally, for all $h, k \in H$ the function $y \to \langle K^\star(s,t,y)K(s,t,y)h,k\rangle$ is constant in $y \in H$.

Invoking Theorem 3.1 and the first part of the proof we obtain for $0 < s < t$ and $y \in H$,

$$\mathbb{E}^y(Z(t)|Z(s)) = m(t,y) + K^\star(s,t,y)Q^{-1/2}(s)(Z(s)-m(s,y)), \quad \mathbb{P}^y - a.s. \tag{3.9}$$

Then the Markov property yields for $0 < s < t$,

$$m(t-s,x) = m(t,y) + K^\star(s,t,y)Q^{-1/2}(s)(x-m(s,y)) \quad \text{for} \quad \mu(s,y)-a.e.\ x, \tag{3.10}$$

Putting $x - m(s,y) = z \in H$ equation (3.10) takes the form

$$m(t-s, z+m(s,y)) = m(t,y) + K^\star(s,t,y)Q^{-1/2}(s)z \quad \text{for} \quad \nu_s - a.e.\ z,$$

where $\nu_s = \mathscr{N}(0, Q(s))$. The measurable linear operator $K^\star(s,t,y)Q^{-1/2}(s)$ is well defined linear on a dense linear space $\operatorname{Im}\left(Q^{1/2}(s)\right)$. By Hypothesis 2.1 the operator

$$K^\star(s,t,y)Q^{-1/2}(s)$$

has a unique extension to a bounded linear operator $L(t,s,y)$ on $H$. Therefore, for any $u, v \in H$ and $x = Q^{1/2}(s)u$ and $z = Q^{1/2}(s)v$ we have

$$m(t-s,x) - m(t-s,z) = K^\star(s,t,y)(u-v).$$

Hence $K^\star(t,s,y)$ is constant in $y \in H$ and (3.4) easily follows. Putting in equation (3.4) $x = m(s,y)$ we obtain (3.5). □

**Lemma 3.3.** *For $0 \leqslant s \leqslant t$ we have $L(t,s) = L(t-s)$ and $(L(t))$ is a strongly continuous semigroup on $H$.*

*Proof.* For any $x, z \in H$ (3.4) yields

$$m(t-s,x) - m(t-s,z) = L(t,s)(x-z). \tag{3.11}$$

Therefore,

$$L(t,s)x = m(t-s,x) - m(t-s,0) \tag{3.12}$$

depends on $t-s$ only and there exists a function, still denoted by $L$, such that

$$L(t,s) = L(t-s), \quad 0 \leqslant s \leqslant t.$$

Invoking (3.11) we obtain

$$\begin{aligned}L(s+t)(x-y) &= m(s+t,x)-m(s+t,y)\\ &= m(s,m(t,x))-m(s,m(t,y)) = L(s)(m(t,x)-m(t,y))\\ &= L(s)L(t)(x-y).\end{aligned}$$

Finally,

$$L(s+t)x = L(s)L(t)x, \quad x \in H,$$

and

$$L(t)x = m(t,x) - m(t,0). \tag{3.13}$$

By Hypothesis 2.1,

$$\lim_{t\downarrow 0}\langle L(t)x, h\rangle = \langle x, h\rangle, \quad h \in H.$$

Since $t \mapsto m(t,y)$ is weakly continuous for every $y \in H$, we obtain

$$\sup_{t\leqslant T}|L(t)y| \leqslant C_T(y).$$

Therefore, the semigroup property of $(L(t))$ and a well known result (see e.g. [18]) imply

$$\lim_{t\downarrow 0} L(t)x = x \quad x \in H,$$

which completes proof of the lemma. □

**Lemma 3.4.** *There exist* $\lambda_0, C > 0$ *such that*

$$|m(t,0)| \leqslant C\mathrm{e}^{\lambda_0 t}, \quad t \geqslant 0.$$

*Proof.* We have

$$m(t,x) = m(t,y) - L(t)y + L(t)x,$$

and therefore, the function

$$g(t) := m(t,y) - L(t)y$$

is independent of $y \in H$. Putting $m(t,0) = g(t)$ it is easy to check that

$$g(t+s) = L(t)g(s) + g(t), \quad s,t \geqslant 0. \tag{3.14}$$

For any $k \geqslant 1$ equation (3.14) yields

$$g(k) = \sum_{i=0}^{k-1} L\left(i\right) g\left(1\right).$$

Since $(L(t))$ is a $C_0$-semigroup in $H$ there exist finite $M$ and $\lambda_0$ such that

$$\|L(t)\| \leqslant M\mathrm{e}^{\lambda_0 t}, \quad t \geqslant 0.$$

Therefore,

$$|g\left(k\right)| \leqslant \left(\sum_{i=0}^{k-1}\|L\left(i\right)\|\right)|g\left(1\right)| \leqslant \left(\sum_{i=0}^{k-1} M\mathrm{e}^{\lambda_0 i}\right)|g\left(1\right)|$$

with
$$c := \frac{M}{\mathrm{e}^{\lambda_0} - 1} |g(1)|.$$
Now, take any $t = k + s$ where $k \geqslant 0$ and $s \in [0, 1)$. Then
$$g(t) = g\left(k + s\right) = L\left(k\right) g(s) + g\left(k\right),$$
and hence
$$|g(t)| \leqslant M\mathrm{e}^{\lambda_0 k} \sup_{t\leqslant 1} |g(s)| + c\mathrm{e}^{\lambda_0 k} \leqslant C\mathrm{e}^{\lambda_0 t},$$
and the lemma follows. □

**Lemma 3.5.** *There exists $b_H \in H$ such that for $b_V = (\lambda - A)b_H \in V$ we have*
$$m(t, 0) = \int_0^t L_V(s) b_V \mathrm{d}s.$$

*Proof.* Let $g(t) = m(t, 0)$. For $\lambda > \lambda_0$ the function
$$\widehat{g}(\lambda) = \int_0^{+\infty} \mathrm{e}^{-\lambda t} g(t) \mathrm{d}t$$
is well defined by Lemma 3.4. Taking the Laplace transform of both sides of equation (3.14) we obtain for $\lambda > \lambda_0$,
$$\mathrm{e}^{\lambda s} \int_0^{\infty} \mathrm{e}^{-\lambda(t+s)} g(t + s) \mathrm{d}t - \widehat{g}(\lambda) = (\lambda - A)^{-1} g(s),$$
hence
$$\mathrm{e}^{\lambda s} \left( \widehat{g}(\lambda) - \int_0^s \mathrm{e}^{-\lambda u} g(u) \mathrm{d}u \right) - \widehat{g}(\lambda) = (\lambda - A)^{-1} g(s).$$
Therefore,
$$\frac{\mathrm{e}^{\lambda s} - 1}{s} \widehat{g}(\lambda) - \frac{\mathrm{e}^{\lambda s}}{s} \int_0^s \mathrm{e}^{-\lambda u} g(u) \mathrm{d}u = (\lambda - A)^{-1} \left( \frac{g(s)}{s} \right).$$
The left-had-side of the above equality converges to
$$b_H = \lambda \widehat{g}(\lambda)$$
for $s \downarrow 0$ and thereby the function $\frac{g(s)}{s}$ has a limit $b_V$ in $V$ for $s \downarrow 0$. More precisely,
$$b_V = \lim_{s \downarrow 0} \frac{g(s)}{s} = (\lambda - A) \left( \lambda \widehat{g}(\lambda) \right).$$
Invoking again equation (3.14) we obtain
$$\lim_{s \downarrow 0} \frac{g(s + t) - g(t)}{s} = \lim_{s \downarrow 0} \left( L_V(t) \frac{g(s)}{s} \right) = L_V(t) b_V,$$
where the convergence holds in $V$. Finally,
$$\frac{\mathrm{d}g}{\mathrm{d}t}(t) = L_V(t) b_V,$$
and the lemma easily follows. □

Let us recall that

$$L(t-s)=K^\star(s,t)Q^{-1/2}(s).$$

Hence Lemma 3.5 implies that (3.3) can be written in the form

$$Q(t+s)=Q(t)+L(t)Q(s)L^\star(t),\quad s,t\geqslant 0. \tag{3.15}$$

It is easy to see that $V'=\mathrm{Dom}\,(A^\star)$ endowed with the graph norm and

$$V'\subset H\subset V.$$

If $R\colon H\to H$ is a trace class operator, such that $R=R^\star\geqslant 0$ then it gives rise to an operator $\tilde R\colon V'\to V$ which is again of trace class, and $\tilde R=\tilde R'\geqslant 0$. In the sequel we will use the same notation $R$ for both operators. For any bounded $R:H\to H$ let

$$\mathscr{L}(t)R=L(t)RL^\star(t).$$

It may be checked that $(\mathscr{L}(t))$ is a $C_0$-semigroup on the space of trace-class operators on $H$. Let $(L_V(t))$ denote an extension of $(L(t))$ to $V$. Then $L_V'(t)$ can be identified with the restriction of $L^\star(t)$ to $V'$. Let $E$ denote a separable Banach space of symmetric trace class operators $R\colon V'\to V$ endowed with the nuclear norm. The dual space $E^\star$ can be identified as the space of bounded operators from $V$ to $V'$ and

$$_{E^\star}\langle P,R\rangle_E=\mathrm{Tr}\,(RP)\,,$$

see pp. 34 and 65 of [8] for details. It is easy to see that the adjoint semigroup acting on $E^\star$ has the form

$$\mathscr{L}^\star(t)R=L'_{-1}(t)RL_V(t),\quad R\colon V\to V'.$$

**Lemma 3.6.** *There exist $\lambda_0>0$ and $C>0$ such that*

$$\|Q(t)\|_E\leqslant C\mathrm{e}^{\lambda_0 t},\quad t\geqslant 0.$$

*Proof.* The proof is similar to the proof of Lemma 3.4 hence omitted. □

**Lemma 3.7.** *There exists $Q\in E$ such that*

$$Q(t)=\int_0^t L_V(s)QL_V'\mathrm{d}s.$$

*Moreover, the restriction of $Q$ to an unbounded operator acting in $H$ (still denoted by $Q$) is selfadjoint in $H$, and for every $t>0$, $L(t)Q^{1/2}$ has an extension to a Hilbert–Schmidt operator on $H$ and*

$$\int_0^t \left\|L(s)Q^{1/2}\right\|_{HS}^2\mathrm{d}s<\infty.$$

*Proof.* Using the same arguments as in the proof of Lemma 3.4 we can show that there exist $\lambda_0, C>0$ such that

$$\|Q(t)\|_E\leqslant C\mathrm{e}^{\lambda_0 t},\quad t\geqslant 0.$$

Therefore, for any $\lambda>\lambda_0$ the operator

$$\widehat{Q}(\lambda)=\int_0^{+\infty}\mathrm{e}^{-\lambda t}Q(t)\mathrm{d}t,$$

where the integral is the Bochner integral in $E$ is a well defined element of $E$. Let $\mathscr{A}$ denote the generator of the semigroup $\mathscr{L}$ on $E$. In view of (3.15) we can apply the same arguments as in the proof of Lemma 3.5 to obtain

$$\lim_{s\downarrow 0}\left((\lambda-\mathscr{A})^{-1}\frac{1}{s}Q(s)\right)=\lambda\widehat{Q}(\lambda), \tag{3.16}$$

where the convergence holds in $E$. Let us recall that for $h,k\in V'$ a bounded operator $h\otimes k\colon V\to V'$ is defined by the formula

$$(h\otimes k)x=\langle h,x\rangle k,\quad x\in V.$$

For $u,v\in \mathrm{Dom}(A)\subset V$ we have $u\otimes v\colon V'\to V$. Moreover, $u\otimes u\in\mathrm{Dom}(\mathscr{A})$ and

$$\mathscr{A}(u\otimes u)=u\otimes(Au)+(Au)\otimes u. \tag{3.17}$$

Similarly, if $u\in V'$ then $u\otimes u\colon V\to V'$, $u\otimes u\in\mathrm{Dom}\left(\mathscr{A}^\star\right)$ and

$$\mathscr{A}^\star(u\otimes u)=u\otimes(A^\star u)+(A^\star u)\otimes u. \tag{3.18}$$

Therefore, for $h\in H\subset V$ and $u\in V'$,

$$\begin{aligned}{}_E\left\langle h\otimes h,\mathscr{A}^\star(u\otimes u)\right\rangle_{E^\star}&={}_E\left\langle h\otimes h,u\otimes(A^\star u)+(A^\star u)\otimes u\right\rangle_{E^\star}\\&={}_E\left\langle h\otimes h,u\otimes(A^\star u)\right\rangle_{E^\star}+{}_E\left\langle h\otimes h,(A^\star u)\otimes u\right\rangle_{E^\star}\\&=\langle h,u\rangle\langle h,A^\star u\rangle+\langle h,A^\star u\rangle\langle h,u\rangle\\&=2\langle h,u\rangle\langle h,A^\star u\rangle\\&=2\left\langle(h\otimes h)u,A^\star u\right\rangle,\end{aligned}$$

where the first step follows from (3.18). Therefore, by the properties of the nuclear norm on $E$ we find that for any $R\in E$ and $u\in V'$,

$${}_E\left\langle R,\mathscr{A}^\star(u\otimes u)\right\rangle_{E^\star}=2\left\langle Ru,A^\star u\right\rangle. \tag{3.19}$$

Equation (3.16) yields

$$\lim_{s\downarrow 0}{}_E\left\langle(\lambda-\mathscr{A})^{-1}\frac{1}{s}Q(s),P\right\rangle_{E^\star}={}_E\left\langle\lambda\widehat{Q}(\lambda),P\right\rangle_{E^\star}$$

for any $P\in E^\star$. Taking $P=(\lambda-\mathscr{A}^\star)(u\otimes u)$ with $u\in\mathrm{Dom}\left(A'_V\right)\subset V'$ and using (3.19) we obtain

$$\begin{aligned}\lim_{s\downarrow 0}\left\langle\frac{1}{s}Q(s)u,u\right\rangle&={}_E\left\langle\lambda\widehat{Q}(\lambda),(\lambda-\mathscr{A}^\star)(u\otimes u)\right\rangle_{E^\star}\\&=\lambda\left\langle\lambda\widehat{Q}(\lambda)u,u\right\rangle-2\lambda\left\langle\widehat{Q}(\lambda)u,A^\star u\right\rangle.\end{aligned}$$

The bilinear form

$$B(u,v)=\lambda\left\langle\lambda\widehat{Q}(\lambda)u,v\right\rangle-\lambda\left\langle\widehat{Q}(\lambda)u,A^\star v\right\rangle-\lambda\left\langle\widehat{Q}(\lambda)v,A^\star u\right\rangle,$$

defined for $u,v\in\mathrm{Dom}\left(A'_V\right)$ is well defined and symmetric as a limit in $s\downarrow 0$ of symmetric bilinear forms. Moreover,

$$|B(u,v)|\leqslant c|u|_{V'}|v|_{V'}.$$

Therefore, there exists $Q\colon V'\to V$ such that

$$_V\left\langle \lim_{s\downarrow 0}\frac{1}{s}Q(s)u,u\right\rangle_{V'} =_V \langle Qu,u\rangle_{V'}.$$

Clearly, $Q$ is non-negative and

$$_V\langle Qu,v\rangle_{V'} =_V \langle Qv,u\rangle_{V'}, \quad u,v\in V'.$$

Therefore, the bilinear form

$$B_H(u,v) =_V \langle Qu,v\rangle_{V'}, \quad u,v\in V'\subset H,$$

defines a selfadjoint operator in $H$, still denoted by $Q$ and such that

$$\mathrm{dom}\left(A^\star\right)\subset \mathrm{dom}\left(Q^{1/2}\right),$$

which proves (2.1). Invoking (3.15) again we readily obtain

$$\frac{\mathrm{d}}{\mathrm{d}t}\langle Q(t)u,u\rangle = 2\left\langle Q^{1/2}A^\star u, Q^{1/2}A^\star u\right\rangle,$$

hence

$$\langle Q(t)u,u\rangle = \int_0^t \left\langle Q^{1/2}L^\star(s)u, Q^{1/2}L^\star(s)u\right\rangle \,\mathrm{d}s, \quad u\in \mathrm{dom}\left(\left(A^\star\right)^2\right).$$

Finally, by polarisation and the density of $\mathrm{dom}\left(\left(A^\star\right)^2\right)$ in $H$ we obtain (2.4). Moreover, $Q\colon V'\to V$ is a trace-class operator since $\widehat{Q}(\lambda)$ is. □

## 4. Proof of Theorem 2.3

We define a $V$-valued and $(\mathscr{F}_t)$-adapted process

$$M(t) = Z(t) - Z(0) - tb_V - \int_0^t A_V Z(s)\,\mathrm{d}s.$$

It is easy to check that

$$M(t) = Z(t) - m(t,x) - \int_0^t A_V(Z(s)-m(s,x))\,\mathrm{d}s, \tag{4.1}$$

where we put $Z(0)=x$. Then $\mathbb{P}^x$-a.s.

$$\begin{aligned}
\mathbb{E}^x\left(M\left(t_2\right)\middle|\,\mathscr{F}_{t_1}\right) &= \mathbb{E}^x\left(Z\left(t_2\right)-m\left(t_2,x\right)\middle|\,Z\left(t_1\right)\right)\\
&\quad - \int_0^{t_1} A_V(Z(s)-m(s,x))\,\mathrm{d}s - \int_{t_1}^{t_2} A_V\mathbb{E}^x\left(Z\left(s\right)-m(s,x)\middle|\,Z\left(t_1\right)\right)\\
&= L\left(t_2-t_1\right)\left(Z\left(t_1\right)-m\left(t_1,x\right)\right)\\
&\quad - \int_0^{t_1} A_V\left(Z(s)-m\left(s,x\right)\right)\,\mathrm{d}s - \int_{t_1}^{t_2} A_V L\left(s-t_1\right)\left(Z\left(t_1\right)-m\left(t_1,x\right)\right)\,\mathrm{d}s\\
&= Z\left(t_1\right)-m\left(t_1,x\right) - \int_0^{t_1} A_V(Z(s)-m(s,x))\,\mathrm{d}s\\
&= M\left(t_1\right).
\end{aligned}$$

Let $M^h(t) = \langle M(t), h\rangle$ for $h \in V' = \mathrm{dom}\,(A^\star)$. Then $M^h$ is a continuous Gaussian martingale with $M^h(0) = 0$. Since the martingale is Gaussian, we have

$$\left\langle M^h\right\rangle_t = \mathbb{E}\left|M^h(t)\right|^2.$$

It remains to compute $\mathbb{E}\left|M^h(t)\right|^2$. To this end we write

$$\begin{aligned}\left|M^h(t)\right|^2 &= \langle Z(t) - m(t,x), h\rangle^2 + \left(\int_0^t \langle Z(s) - m(s,x), A^\star h\rangle\, \mathrm{d}s\right)^2\\ &\quad - 2\int_0^t \langle Z(t) - m(t,x), h\rangle \langle Z(s) - m(s,x), A^\star h\rangle\, \mathrm{d}s.\end{aligned}$$

Then

$$\begin{aligned}\mathbb{E}\left|M^h(t)\right|^2 &= \langle Q(t)h, h\rangle + \int_0^t\int_0^t \mathbb{E}\,\langle Z(s) - m(s,x), A^\star h\rangle \langle Z(u) - m(u,x), A^\star h\rangle\, \mathrm{d}u\mathrm{d}s\\ &\quad - 2\int_0^t \langle L(t-s)Q(s)A^\star h, h\rangle\, \mathrm{d}s.\end{aligned}$$

Since $\langle Q(t)h, h\rangle$ can be written in the form

$$\langle Q(t)h, h\rangle = \int_0^t \left\langle Q^{1/2}L^\star(t-s)h, Q^{1/2}L^\star(t-s)h\right\rangle \mathrm{d}s$$

we obtain

$$\begin{aligned}\frac{d}{dt}\mathbb{E}\left|M^h(t)\right|^2 &= 2\,\langle Q(t)h, A^\star h\rangle + \left|Q^{1/2}h\right|^2\\ &\quad + \int_0^t \mathbb{E}\,\langle Z(s) - m(s,x), A^\star h\rangle \langle Z(t) - m(t,x), A^\star h\rangle\, \mathrm{d}s\\ &\quad + \int_0^t \mathbb{E}\,\langle Z(t) - m(t,x), A^\star h\rangle \langle Z(u) - m(u,x), A^\star h\rangle\, \mathrm{d}u\\ &\quad - 2\int_0^t \langle L(t-s)Q(s)A^\star h, A^\star h\rangle\, \mathrm{d}s - 2\,\langle Q(t)h, A^\star h\rangle\\ &= \left|Q^{1/2}h\right|^2.\end{aligned}$$

Finally, by the martingale representation theorem, see e.g. [6], there exists a standard cylindrical Wiener process $W$ on $H$ such that $M^h(t) = \left\langle W(t), Q^{1/2}h\right\rangle$ is adapted to $(\mathscr{F}_t)$ and for any $h \in \mathrm{dom}\,(A^\star)$,

$$\langle Z(t), h\rangle = \langle Z(0), h\rangle + \int_0^t \langle Z(s), A^\star h\rangle\, \mathrm{d}s + \int_0^t \langle b_V, h\rangle\, \mathrm{d}s + \left\langle W(t), Q^{1/2}h\right\rangle.$$

Now, Theorem 2.3 follows from a result in [4].

## 5. Proof of Theorem 2.4

It was proved in [17] that $P_tC_b\left(H_{bw}\right)\subset C_b\left(H_{bw}\right)$. Let $B_r\subset H$ denote the centered closed ball of radius $r$. Let us recall that bounded sets of $H_{bw}$ are precisely the bounded sets of the norm topology and the sets $B_r$ are compact in $H_{bw}$. The proof is a simple modification of the proof of Theorem 4.2 in [12]. Let us note here that by definition of the $bw$-topology, see [10] the space $H_{bw}$ is a $k$-space and therefore the space $C_b^{strict}\left(H_{bw}\right)$ is complete. We provide details only for the proof of strong continuity. To this end it is enough to show that for every weakly sequentially continuous $\phi : H\to\mathbb{R}$ and every $r>0$,

$$\lim_{t\to 0}\sup_{x\in B_r}\left|P_t\phi(x)-\phi(x)\right|=0.$$

Assume that this convergence does not hold. Then there exist $\varepsilon>0$, a sequence $t_n\to 0$ and a sequence $\left(x_n\right)\subset B_r$ such that

$$\left|P_{t_n}\phi\left(x_n\right)-\phi\left(x_n\right)\right|>\varepsilon.$$

Therefore

$$\int_H\left|\phi\left(m\left(t_n,x_n\right)+y\right)-\phi\left(x_n\right)\right|\mu\left(t_n,dy\right)>\varepsilon, \tag{5.1}$$

where by Lemma 2.2 the measure $\mu(t,dy)=\mathscr{N}(0,Q(t))(dy)$ is independent of $x\in H$. Hypothesis 2.1 yields

$$\left\langle m\left(t_n,0\right),h\right\rangle\to 0,\quad h\in H,$$

and clearly

$$\lim_{n\to\infty}\left\langle x_n,L^\star\left(t_n\right)h\right\rangle=\langle x,h\rangle,\quad h\in H.$$

Note that (2.3) yields

$$\langle m(t,0),h\rangle=\left\langle(\lambda-A)\int_0^t L(s)b_H\,\mathrm{d}s,h\right\rangle$$

and therefore (2.3) holds for all $h\in H$. Finally,

$$\left\langle m\left(t_n,x_n\right),h\right\rangle=\left\langle x_n,L^\star\left(t_n\right)h\right\rangle+\left\langle m\left(t_n,0\right),h\right\rangle\longrightarrow\langle x,h\rangle,\quad h\in H.$$

Now, using the dominated convergence, we can pass to the limit in (5.1) obtaining the desired contradiction.

We pass now to the proof of the second part of the theorem. It is easy to check that for any $\phi\in\mathscr{F}C_b^\infty\left(A^\star\right)$ formula (5.1) holds and by the definition of the space $\mathscr{F}C_b^\infty\left(A^\star\right)$ the function $L\phi$ is a well defined element of $C_b^{strict}\left(H_{bw}\right)$. The proof that $\mathscr{F}C_b^\infty\left(A^\star\right)$ is a core for $L$ is an easy modification of the proofs of Lemma 4.4 and Theorem 4.5 in [12] and thus omitted.

## 6. Example: boundary noise

We will consider here a stochastic PDE with boundary noise introduced in [1], see also [11]:

$$\begin{cases} \frac{\partial u}{\partial t}(t,\xi) = \frac{\partial^2 u}{\partial \xi^2}(t,\xi), & t>0, \xi>0, \\ u(t,0) = \dot{W}(t), & t>0, \\ u(0,\xi) = x(\xi), & \xi \geqslant 0. \end{cases} \tag{6.1}$$

In [1] solution to (6.1) is defined for all $t>0$ and $\xi \geqslant 0$ by the formula

$$u(t,\xi) = \int_0^{+\infty} G(t,\xi,\eta)x(\eta)\,\mathrm{d}\eta + \int_0^t \frac{\partial G}{\partial \eta}(t-s,\xi,0)\,\mathrm{d}W(s), \tag{6.2}$$

where

$$G(t,\xi,\eta) = \frac{1}{\sqrt{4\pi t}}\left(\mathrm{e}^{-|\xi-\eta|^2/4t} - \mathrm{e}^{-|\xi+\eta|^2/4t}\right), \quad \xi,\eta \geqslant 0.$$

Let $\rho(\xi) = \min\left(1, \xi^{1+\alpha}\right)$ with $\alpha \in (0,1)$ and $H = L^2([0,+\infty), \rho(\xi)\mathrm{d}\xi)$. Finally, let $A$ denote the Dirichlet Laplacian in $L^2([0,+\infty))$.

**Proposition 6.1.** *(1) The operator $A$ has an extension to the generator (still denoted by $A$) of an analytic $C_0$-semigroup $\left(\mathrm{e}^{tA}\right)$ on $H$.*
*(2) There exists a selfadjoint unbounded operator $Q$ in $H$, such that*

$$u(t) = \mathrm{e}^{tA}x + \int_0^t \mathrm{e}^{(t-s)A}Q^{1/2}\,\mathrm{d}W(s).$$

*Proof.* Part 1 of the proposition has been proved in [15]. The results in [1] imply that for any $T>0$,

$$\sup_{t\leqslant T} \mathbb{E}|u(t)|_H^2 < \infty, \tag{6.3}$$

and obviously the process is Gaussian. Clearly,

$$m(t,x) = \mathrm{e}^{tA}x, \quad t \geqslant 0,\, x \in H,$$

and

$$Q(t)x(\xi) = \int_0^\infty q(t,\xi,\eta)\rho(\eta)x(\eta)\,\mathrm{d}\eta, \quad x \in H,\, \xi \in [0,+\infty),$$

where

$$q(t,\xi,\eta) = \int_0^t \frac{\partial G}{\partial \eta}(s,\xi,0)\frac{\partial G}{\partial \eta}(s,\eta,0)\,\mathrm{d}s.$$

It is not very difficult to show that $\mathbb{E}\,|u(t)|_H^2 < \infty$, and that $u$ is a mean-square continuous, Gaussian and Markov process in $H$. For a fixed $\lambda > 0$ we denote by $D$ the Dirichlet map, e.a. for $a > 0$ we define $Da$ as a solution to the equation $(\lambda - A)\,\phi = 0$, that is

$$\left(\lambda - \frac{\partial^2}{\partial \xi^2}\right)\phi = 0, \quad \phi(0) = a.$$

It is easy to see that

$$Da(\xi) = a\mathrm{e}^{-\xi\sqrt{\lambda}}, \quad \xi \geqslant 0.$$

Let $H^{-s}(0,+\infty)$ denote the negative Sobolev space. Using arguments similar to those in [7] one can show that a unique mild solution to equation (6.1) exists in $H^{-s}(0,+\infty)$ provided $s > \frac{1}{2}$. Moreover,

$$u(t) = \mathrm{e}^{tA}x + (\lambda - A)\int_0^t \mathrm{e}^{(t-s)A} D\mathrm{d}W(s).$$

□

School of Mathematics and Statistics, The University of Sydney, Sydney 2006, Australia

*E-mail address*: Beniamin.Goldys@sydney.edu.au

Institute of Mathematics, Polish Academy of Sciences, Św. Tomasza 30/7, 31-027 Cracow, Poland

*E-mail address*: napeszat@cyf-kr.edu.pl

Institute of Mathematics, Polish Academy of Sciences, Śniadeckich 8, 00-950 Warsaw, Poland

*E-mail address*: J.Zabczyk@impan.pl